\documentclass[twoside,leqno,10pt, A4]{amsart}
\usepackage{amsfonts}
\usepackage{amsmath}
\usepackage{amscd}
\usepackage{amssymb}
\usepackage{amsthm}
\usepackage{amsrefs}
\usepackage{latexsym}
\usepackage{mathrsfs}
\usepackage{bbm}
\usepackage{enumerate}
\usepackage{graphicx}

\usepackage{amsfonts}
\usepackage{amsmath}
\usepackage{amscd}
\usepackage{amssymb}
\usepackage{amsthm}
\usepackage{amsrefs}
\usepackage{latexsym}
\usepackage{mathrsfs}
\usepackage{bbm}
\usepackage{amscd}
\usepackage{amssymb}
\usepackage{amsthm}
\usepackage{amsrefs}
\usepackage{latexsym}
\usepackage{mathrsfs}
\usepackage{bbm}
\usepackage{enumerate}
\usepackage{graphicx}
\usepackage{color}
\begin{document}

\newtheorem{theorem}[subsection]{Theorem}
\newtheorem{proposition}[subsection]{Proposition}
\newtheorem{lemma}[subsection]{Lemma}
\newtheorem{corollary}[subsection]{Corollary}
\newtheorem{conjecture}[subsection]{Conjecture}
\newtheorem{prop}[subsection]{Proposition}
\numberwithin{equation}{section}
\newcommand{\mr}{\ensuremath{\mathbb R}}
\newcommand{\mc}{\ensuremath{\mathbb C}}
\newcommand{\dif}{\mathrm{d}}
\newcommand{\intz}{\mathbb{Z}}
\newcommand{\ratq}{\mathbb{Q}}
\newcommand{\natn}{\mathbb{N}}
\newcommand{\comc}{\mathbb{C}}
\newcommand{\rear}{\mathbb{R}}
\newcommand{\prip}{\mathbb{P}}
\newcommand{\uph}{\mathbb{H}}
\newcommand{\fief}{\mathbb{F}}
\newcommand{\majorarc}{\mathfrak{M}}
\newcommand{\minorarc}{\mathfrak{m}}
\newcommand{\sings}{\mathfrak{S}}
\newcommand{\fA}{\ensuremath{\mathfrak A}}
\newcommand{\mn}{\ensuremath{\mathbb N}}
\newcommand{\mq}{\ensuremath{\mathbb Q}}
\newcommand{\half}{\tfrac{1}{2}}
\newcommand{\f}{f\times \chi}
\newcommand{\summ}{\mathop{{\sum}^{\star}}}
\newcommand{\chiq}{\chi \bmod q}
\newcommand{\chidb}{\chi \bmod db}
\newcommand{\chid}{\chi \bmod d}
\newcommand{\sym}{\text{sym}^2}
\newcommand{\hhalf}{\tfrac{1}{2}}
\newcommand{\sumstar}{\sideset{}{^*}\sum}
\newcommand{\sumprime}{\sideset{}{'}\sum}
\newcommand{\sumprimeprime}{\sideset{}{''}\sum}
\newcommand{\sumflat}{\sideset{}{^\flat}\sum}
\newcommand{\shortmod}{\ensuremath{\negthickspace \negthickspace \negthickspace \pmod}}
\newcommand{\V}{V\left(\frac{nm}{q^2}\right)}
\newcommand{\sumi}{\mathop{{\sum}^{\dagger}}}
\newcommand{\mz}{\ensuremath{\mathbb Z}}
\newcommand{\leg}[2]{\left(\frac{#1}{#2}\right)}
\newcommand{\muK}{\mu_{\omega}}
\newcommand{\thalf}{\tfrac12}
\newcommand{\lp}{\left(}
\newcommand{\rp}{\right)}
\newcommand{\Lam}{\Lambda_{[i]}}
\newcommand{\lam}{\lambda}
\def\L{\fracwithdelims}
\def\om{\omega}
\def\pbar{\overline{\psi}}
\def\phis{\phi^*}
\def\lam{\lambda}
\def\lbar{\overline{\lambda}}
\newcommand\Sum{\Cal S}
\def\Lam{\Lambda}
\newcommand{\sumtt}{\underset{(d,2)=1}{{\sum}^*}}
\newcommand{\sumt}{\underset{(d,2)=1}{\sum \nolimits^{*}} \widetilde w\left( \frac dX \right) }

\newcommand{\hf}{\tfrac{1}{2}}
\newcommand{\af}{\mathfrak{a}}
\newcommand{\Wf}{\mathcal{W}}

\theoremstyle{plain}
\newtheorem{conj}{Conjecture}
\newtheorem{remark}[subsection]{Remark}

\makeatletter
\def\widebreve{\mathpalette\wide@breve}
\def\wide@breve#1#2{\sbox\z@{$#1#2$}%
     \mathop{\vbox{\m@th\ialign{##\crcr
\kern0.08em\brevefill#1{0.8\wd\z@}\crcr\noalign{\nointerlineskip}%
                    $\hss#1#2\hss$\crcr}}}\limits}
\def\brevefill#1#2{$\m@th\sbox\tw@{$#1($}%
  \hss\resizebox{#2}{\wd\tw@}{\rotatebox[origin=c]{90}{\upshape(}}\hss$}
\makeatletter

\title[Lower Bounds for negative moments of Dirichlet $L$-functions to a fixed modulus]{Lower bounds for negative moments of Dirichlet $L$-functions to a fixed modulus}

\author[P. Gao]{Peng Gao}
\address{School of Mathematical Sciences, Beihang University, Beijing 100191, China}
\email{penggao@buaa.edu.cn}

\begin{abstract}
 We establish lower bounds for the $2k$-th moment of central values of the family of primitive Dirichlet $L$-functions to a fixed prime modulus for all real $k<0$, assuming the non-vanishing of these $L$-values.

\end{abstract}

\maketitle



\section{Introduction}
\label{sec 1}

 The generalized Riemann hypothesis (GRH) asserts that all non-trivial zeros of Dirichlet $L$-functions can be written as $\rho = 1/2+i\gamma$ with $\gamma \in \mr$. Moreover, it is believed that there are no $\mq$-linear relations among the non-negative $\gamma$'s. In particular, this implies that $L(\half, \chi) \neq 0$ for all primitive Dirichlet characters $\chi$ and it is known as a conjecture of S. Chowla \cite{chow} when $\chi$ is quadratic.

   One way to investigate the non-vanishing problem is to evaluate the one-level densities of low-lying zeros of families of $L$-functions. In fact, the density conjecture of N. Katz and P. Sarnak \cites{KS1, K&S} implies that $L(1/2, \chi) \neq 0$ for almost all primitive Dirichlet $L$-functions.  Computing essentially the one-level density of low-lying zeros of the corresponding families of Dirichlet $L$-functions for test functions whose Fourier transforms being supported in $(-2, 2)$, M. R. Murty  \cite{Murty} showed that under GRH, at least $50\%$ of the primitive Dirichlet $L$-functions do not vanish at the central point. See also the work of H. P. Hughes and Z. Rudnick In \cite{HuRu} regarding the one-level density of low-lying zeros of the family of primitive Dirichlet $L$-functions to a fixed prime modulus.

  Another way to address whether $L(\half, \chi)=0$ or not is to study moments of families of $L$-functions. In \cite{BM}, B. Balasubramanian and V. K. Murty showed that $L(1/2,\chi) \neq 0$ for at least $4\%$ of Dirichlet characters $\chi$ to a fixed modulus $q$ by evaluating the first and second mollified moments of $L(1/2, \chi)$. For a prime $q$, this proportion was improved to $1/3$ by H. Iwaniec and P. Sarnak \cite{I&S}, to $34.11\%$ by H. M. Bui \cite{Bui}, to $3/8$ by R. Khan and H. T. Ngo \cite{KN}, and to $5/13$ by R. Khan, D. Mili\'cevi\'c and H. T. Ngo \cite{KMN22}. 

 Because of the important role played by moments of $L$-functions concerning the non-vanishing issue, a considerable amount of work has been done in this direction. For the family of Dirichlet $L$-functions to a fixed modulus $q \not \equiv 2 \pmod 4$ (to ensure primitive Dirichlet
  characters modulo $q$ exist), it is widely believed that (see \cite{R&Sound})  for all real $k  \geq 0$,
\begin{align}
\label{moments}
 \sumstar_{\substack{ \chi \shortmod q }}|L(\tfrac{1}{2},\chi)|^{2k} \sim C_k \phis(q)(\log q)^{k^2},
\end{align}
  where the numbers $C_k$ are explicit constants, $\phis(q)$ denotes the number of primitive characters modulo $q$ and where we denote $\sumstar$ the sum over primitive Dirichlet characters modulo $q$ throughout the paper.

 The formula given in \eqref{moments} was conjectured by J. B. Conrey, D. W. Farmer, J. P. Keating,
  M. O. Rubinstein and N. C. Snaith in \cite{CFKRS} for all positive integral values of $k$ and is well-known for $k=1$. The case $k=2$ was established by D. R. Heath-Brown \cite{HB81} for almost all $q$ and was later shown to hold for all $q$ by K. Soundararajan \cite{Sound2007}.  Subsequent improvements on the error terms in Soundararajan's result can be found in \cites{Young2011, BFKMM1, BFKMM, Wu2020, BPRZ}.

  Although it is challenging to prove \eqref{moments} even for integers $k \geq 3$, much progress has been made for building upper and lower bounds of the conjectured order of magnitude
   for the expression on the left-hand side of \eqref{moments}. In fact, there are now several systematic approaches towards establishing sharp lower and upper bounds. Notably, there are the upper bounds principle due to  M. Radziwi{\l\l} and K. Soundararajan  \cite{Radziwill&Sound} as well as the  lower bounds principle due to W. Heap and K. Soundararajan  \cite{H&Sound}. Other methods can be found in \cites{Sound01,HB2010, Harper, R&Sound1, Radziwill&Sound1, C&L, Gao2024-6}. These results together imply that for large prime $q$ and any real number $k \geq 0$, we have
\begin{align}
\label{orderofmag}
   \sumstar_{\substack{ \chi \shortmod q }}|L(\tfrac{1}{2},\chi)|^{2k} \asymp \phis(q)(\log q)^{k^2}.
\end{align}
  Here we note that one needs to assume GRH in order to show that $\displaystyle \sumstar_{\substack{ \chi \pmod q }}|L(\tfrac{1}{2},\chi)|^{2k} \ll_k \phis(q)(\log q)^{k^2}$ for $k>1$.

  As the order of magnitude for the non-negative moments of the family of Dirichlet $L$-functions to a fixed prime modulus are established in  \eqref{orderofmag}, it is natural to turn to the negative moments of this family and we assume that $L(\tfrac{1}{2},\chi) \neq 0$ for each $\chi$ in the family in order for negative powers of $|L(\tfrac{1}{2},\chi)|$ to be meaningful.  An analogue case for the family of quadratic Dirichlet $L$-functions has been considered by the author in \cite{Gao2022-1}. As already being pointed out in \cite{Gao2022-1},  the behaviour of the negative moments may be more difficult to predict compared to that of the positive ones. One can also see this from the case $\beta=2$ of \cite[Corollary 1]{FK}, where computations by P. J. Forrester and J. P. Keating based on random matrix theory suggests certain phase changes in the asymptotic formulas for the $2k$-th moment of the family of Dirichlet $L$–functions to a fixed modulus when $2k = -(2j- 1)$ for any positive integer $j$.

 In this paper, we establish lower bounds for negative moments of the family of primitive Dirichlet $L$-functions to a fixed prime modulus. Our result is as follows.
\begin{theorem}
\label{thmlowerbound}
   Let $q$ be a large prime number and assume that $L(\half, \chi) \neq 0$ for any primitive Dirichlet character $\chi \pmod q$. Then we have for any real number $k < 0$,
\begin{align*}
   \sumstar_{\substack{ \chi \shortmod q }}|L(\tfrac{1}{2},\chi)|^{2k} \gg_k \phis(q)(\log q)^{k^2}.
\end{align*}
\end{theorem}

  We note that if one interprets $0^{2k}=+\infty$ for $k<0$, then the statement of Theorem \ref{thmlowerbound} is still valid without the assumption that $L(\half, \chi) \neq 0$.  
  We shall derive Theorem \ref{thmlowerbound} by applying a variant of the lower bounds principle of W. Heap and K. Soundararajan \cite{H&Sound}. Such an approach has already been employed by W. Heap, J. Li and J. Zhao \cite{HLZ}, as well as by P. Gao and L. Zhao \cite{G&Zhao2023} to study lower bounds of discrete negative moments of the derivative of the Riemann zeta function $\zeta(s)$ at non-trivial zeros. We expect that the bounds given in Theorem \ref{thmlowerbound} are sharp for $-1 < 2k <0$, by taking account of the prediction in \cite{FK} based on random matrix theory .

\section{Proof of Theorem \ref{thmlowerbound}}
\label{sec 2'}

\subsection{Setup}

  For any real number $\ell$, we let $\lceil \ell \rceil = \min \{ m \in \intz : \ell \leq m \}$ be the celing function of $\ell$. We define a sequence of even natural numbers $\{ \ell_j \}_{1 \leq j \leq R}$ such that $\ell_1= 2\lceil N \log \log q\rceil$ and $\ell_{j+1} = 2 \lceil N \log \ell_j \rceil$ for $j \geq 1$, where $N, M$ are two large natural numbers depending on $k$ only, and where $R$ is
  defined to be the largest natural number satisfying $\ell_R >10^M$. We choose $M$ large enough to ensure that $\ell_{j} >
  \ell_{j+1}^2$ for all $1 \leq j \leq R-1$. It follows from this that 
\begin{align}
\label{sumoverell}
  \sum^R_{j=1}\frac 1{\ell_j} \leq \frac 2{\ell_R}.
\end{align}

 Let ${P}_1$ be the set of odd primes not exceeding $q^{1/\ell_1^2}$ and let
${ P_j}$ be the set of primes lying in the interval $(q^{1/\ell_{j-1}^2}, q^{1/\ell_j^2}]$ for $2\le j\le R$. We
write for each $1 \leq j \leq R$, 
\begin{equation*}
{\mathcal P}_j(\chi) = \sum_{p\in P_j} \frac{1}{\sqrt{p}} \chi(p), \quad  {\mathcal Q}_j(\chi, k) =\Big (\frac{12 (1+|k|) {\mathcal
P}_j(\chi) }{\ell_j}\Big)^{(2-k/(1-k))\ell_j}.
\end{equation*}
  We further define ${\mathcal Q}_{R+1}(\chi, k)=1$.

  We define for any non-negative integer $\ell$ and any complex number $x$,
\begin{equation*}
E_{\ell}(x) = \sum_{j=0}^{\ell} \frac{x^{j}}{j!}.
\end{equation*}
  Further, we define for each $1 \leq j \leq R$ and any real number $\alpha$,
\begin{align}
\label{defN}
{\mathcal N}_j(\chi, \alpha) = E_{\ell_j} (\alpha {\mathcal P}_j(\chi)), \quad \mathcal{N}(\chi, \alpha) = \prod_{j=1}^{R} {\mathcal
N}_j(\chi,\alpha).
\end{align}
   
   Similarly, we define for each $1 \leq j \leq R$ and any real number $\alpha$,
\begin{align}
\label{defM}
{\mathcal M}_j(\chi, \alpha) = E_{\ell_j} (\alpha \Re{\mathcal P}_j(\chi)), \quad \mathcal{M}(\chi, \alpha) = \prod_{j=1}^{R} {\mathcal
M}_j(\chi,\alpha).
\end{align}   

  We now present a variant in our setting of the lower bounds principle of W. Heap and K. Soundararajan \cite{H&Sound}.  To do so, we first note that as $\ell_j$ is even for each $j$, it follows from \cite[Lemma 1]{Radziwill&Sound} that ${\mathcal M}_j(\chi, \alpha)>0$ for any real $\alpha$, which then implies that $\mathcal{M}(\chi, \alpha)>0$ as well. Moreover, it follows from \cite[Lemma 4.1]{Gao2021-2} that for any real number $\alpha$,
\begin{align*}
 \mathcal{M}(\chi, \alpha)\mathcal{M}(\chi, -\alpha)  \geq 1.
\end{align*}

  We apply the above to see for any $c>0$, 
\begin{align}
\label{firstdecomp}
\begin{split}
 \sumstar_{\substack{ \chi \shortmod q }}\Big |\mathcal{N}(\chi, k)\Big |^2 
 \leq & \sumstar_{\substack{ \chi \shortmod q }}\Big |\mathcal{N}(\chi, k)\Big |^2 \Big (\mathcal{M}(\chi, k-1)\mathcal{M}(\chi, 1-k)\Big )^c\\
 =& \sumstar_{\substack{ \chi \shortmod q }}L(\frac 12, \chi)^{-c}\cdot \big(L(\frac 12, \chi)\mathcal{M}(\chi, k-1)\big)^{c} \cdot \Big |\mathcal{N}(\chi, k)\Big |^2 \mathcal{M}(\chi, 1-k)^c.
\end{split}
\end{align}

As $k<0$,  we fix a constant $0<c<-2k/(1-k)$ satisfying
\begin{align}
\label{ccond}
\begin{split}
& 0< \frac {c}{2}-\frac {c}{2k} <1. 
\end{split}
\end{align} 

  We then apply H\"older's inequality with exponents $-2k/c, 2/c, (1+(1-k)c/(2k))^{-1}$ to the right-hand side expression in \eqref{firstdecomp} to see that
\begin{align}
\label{basiclowerbound}
\begin{split}
\sumstar_{\substack{ \chi \shortmod q }}\Big |\mathcal{N}(\chi, k)\Big |^2 \leq & \Big ( \sumstar_{\substack{ \chi \shortmod q }}|L(\tfrac{1}{2},\chi)|^{2k} \Big )^{-c/(2k)}\Big ( \sumstar_{\substack{ \chi \shortmod q
 }}|L(\tfrac{1}{2},\chi)|^2 |\mathcal{M}(\chi, k-1)|^2  \Big)^{c/2} \\
 & \times \Big ( \sumstar_{\substack{ \chi \shortmod q }}   \Big |\mathcal{N}(\chi,
 k)^2\mathcal{M}(\chi, 1-k)^c\Big |^{(1+(1-k)c/(2k))^{-1}}
 \Big)^{1+(1-k)c/(2k)}.
\end{split}
\end{align}
 
   In particular, we set $c=-k/(1-k)$ to see that the condition \eqref{ccond} is satisfied. It then follows from \eqref{basiclowerbound} that we have
\begin{align}
\label{basiclowerbound1}
\begin{split}
 \sumstar_{\substack{ \chi \shortmod q }}\Big |\mathcal{N}(\chi, k)\Big |^2
 \leq & \Big ( \sumstar_{\substack{ \chi \shortmod q }}|L(\tfrac{1}{2},\chi)|^{2k} \Big )^{1/(2(1-k))}\Big ( \sumstar_{\substack{ \chi \shortmod q
 }}|L(\tfrac{1}{2},\chi)|^2 |\mathcal{M}(\chi, k-1)|^2  \Big)^{-k/(2(1-k))} \\
 & \times \Big ( \sumstar_{\substack{ \chi \shortmod q }}   \Big |\mathcal{N}(\chi,
 k)^2\mathcal{M}(\chi, 1-k)^{-k/(1-k)}\Big |^{2}
 \Big)^{1/2}.
\end{split}
\end{align}

   We deduce from \eqref{basiclowerbound1} that in order to prove Theorem \ref{thmlowerbound},
it suffices to establish the following three propositions.
\begin{proposition}
\label{Prop4} With the notation as above, we have
\begin{align*}
\sumstar_{\substack{ \chi \shortmod q }}\big|\mathcal{N}(\chi, k)  \big |^2 \gg_k \phis(q)(\log q)^{ k^2
} .
\end{align*}
\end{proposition}

\begin{proposition}
\label{Prop5} With the notation as above, we have
\begin{align*}
\sumstar_{\substack{ \chi \shortmod q }}   \Big |\mathcal{N}(\chi,
 k)^2\mathcal{M}(\chi, 1-k)^{-k/(1-k)}\Big |^{2}  \ll_k
\phis(q)(\log q)^{ k^2 }.
\end{align*}
\end{proposition}

\begin{proposition}
\label{Prop6} With the notation as above, we have
\begin{align*}
 \sumstar_{\substack{ \chi \shortmod q }} \big|L(\half, \chi){\mathcal M}(\chi, k-1) \big |^{2}  \ll_k
\phis(q)(\log q)^{ k^2 }.
\end{align*}
\end{proposition}

  We shall prove the above Propositions in the rest of the paper.  

\subsection{Proof of Proposition \ref{Prop4}}

    We denote $w(n)$ the multiplicative function satisfying
    $w(p^{\alpha}) = \alpha!$ for prime powers $p^{\alpha}$.  We also define functions $b_j(n), 1 \leq j \leq R$ such that $b_j(n)=1$ when $\Omega(n) \leq \ell_j$ and the primes dividing $n$ are all from the interval $P_j$, where $\Omega(n)$ denotes the number of prime powers dividing $n$. For other values of $n$, we set $b_j(n)=0$. Using these notations, we see that for
    any real number $\alpha$,
\begin{equation}
\label{5.1}
{\mathcal N}_j(\chi, \alpha) = \sum_{n_j} \frac{1}{\sqrt{n_j}} \frac{\alpha^{\Omega(n_j)}}{w(n_j)}  b_j(n_j) \chi(n_j), \quad 1\le j\le R.
\end{equation}
    Observe that each ${\mathcal N}_j(\chi, \alpha)$ is a short Dirichlet polynomial as $b_j(n_j)=0$ unless $n_j \leq
    (q^{1/\ell_j^2})^{\ell_j}=q^{1/\ell_j}$. This together with \eqref{sumoverell} implies that ${\mathcal N}(\chi, k)$ is also a short Dirichlet
    polynomial whose lengths does not exceed $q^{1/\ell_1+ \ldots +1/\ell_R} < q^{2/10^{M}}$. Also, we see from \eqref{5.1}
  that for each $\chi$ modulo $q$, including the case with $\chi$ being the principal character $\chi_0$ modulo $q$,
\begin{align*}
 & |{\mathcal N}(\chi, k)|^2 \ll q^{2(1/\ell_1+ \ldots +1/\ell_R)} < q^{4/10^{M}}.
\end{align*}

   It follows from this that we have
\begin{align}
\label{upperboundprodofN}
\begin{split}
 & \sumstar_{\substack{ \chi \shortmod q }}|{\mathcal N}(\chi, k)|^2 \geq  \sum_{\substack{ \chi \shortmod q }}|{\mathcal N}(\chi, k)|^2+O(q^{4/10^{M}}).
\end{split}
\end{align}

 Denote $\phi(q)$ the  Euler totient function, we recall that the orthogonality relation for Dirichlet characters (see \cite[Corollary 4.5]{MVa1}) asserts 
\begin{align}
\label{orthrel}
 \sum_{\substack{ \chi \shortmod q }} \chi(n)=\begin{cases}
 \varphi(q) \quad \text{if} \ n \equiv 1 \pmod q, \\
 0 \quad \text{otherwise}.
\end{cases}
\end{align}

  As the length of the Dirichlet series ${\mathcal N}(\chi, k)$ is $\leq q$, we deduce from \eqref{orthrel} that only the
  diagonal terms in the last sum of \eqref{upperboundprodofN} survive. Thus, we have
\begin{align*}
\begin{split}
 & \sum_{\substack{ \chi \shortmod q }}|{\mathcal N}(\chi, k)|^2 \ge  \phi(q) \prod^R_{j=1} \Big (  \sum_{n_j} \frac{k^{2\Omega(n_j)}}{n_j w^2(n_j)}  b_j(n_j)\Big )\geq \phis(q) \prod^R_{j=1} \Big (  \sum_{n_j} \frac{k^{2\Omega(n_j)}}{n_j w^2(n_j)}  b_j(n_j)\Big ) \gg   \phis(q)(\log q)^{k^2},
\end{split}
\end{align*}
  where the last estimation above follows from \cite[(6.3)]{Gao2024-6} and the treatments given in \cite[Section 6]{Gao2024-6}. This together with \eqref{upperboundprodofN} now implies the assertions of Proposition \ref{Prop4}.

\subsection{Proof of Proposition \ref{Prop5}}

  We notice that it follows from  \cite[(3.5)]{Gao2024-6} that we have for $|z| \le aK/10$ with $0<a \leq 1$,
\begin{align}
\label{Ebound}
\Big| \sum_{r=0}^K \frac{z^r}{r!} - e^z \Big| \le \frac{|z|^{K}}{K!} \le \Big(\frac{a e}{10}\Big)^{K}.
\end{align}
   
  For any fixed $1 \leq j \leq R$, we apply \eqref{Ebound} with $z=k{\mathcal P}_j(\chi), K=\ell_j$ and $a=1$ to see that when $|{\mathcal P}_j(\chi)| \le \ell_j/(10(1+|k|))$,
\begin{align}
\label{Njboundk}
{\mathcal N}_j(\chi, k)=& \exp( k {\mathcal P}_j(\chi))\Big( 1+  O\Big(\exp(  |k{\mathcal P}_j(\chi)|)\Big(\frac{e}{10}\Big)^{\ell_j} \Big )
=  \exp( k {\mathcal P}_j(\chi))\Big( 1+  O\Big( e^{-\ell_j} \Big )\Big ).
\end{align}
  Similarly, when $|{\mathcal P}_j(\chi)| \le \ell_j/10(1+|k|)$, we have
\begin{align}
\label{Njboundk1}
{\mathcal M}_j(\chi, 1-k)= & \exp( (1-k) \Re{\mathcal P}_j(\chi))\Big( 1+  O\Big(e^{-\ell_j} \Big ) \Big ).
\end{align}

  It follows from \eqref{Njboundk} and \eqref{Njboundk1} that when $|{\mathcal P}_j(\chi)| \le \ell_j/(10(1+|k|))$,
\begin{align}
\label{est1}
 \Big |\mathcal{N}_j(\chi,
 k)^2\mathcal{M}_j(\chi, 1-k)^{-k/(1-k)}\Big |^{2} 
=& \exp( 2k \Re ({\mathcal P}_j(\chi)))\Big( 1+ O\big( e^{-\ell_j} \big) \Big) = |{\mathcal N}_j(\chi, k)|^2 \Big( 1+ O\big(e^{-\ell_j} \big)
\Big).
\end{align}

  On the other hand, we notice that when $|{\mathcal P}_j(\chi)| \ge \ell_j/(10(1+|k|))$,
\begin{align}
\label{4.2}
\begin{split}
|{\mathcal N}_j(\chi, k)| &\le \sum_{r=0}^{\ell_j} \frac{|k{\mathcal P}_j(\chi)|^r}{r!} \le
|{\mathcal P}_j(\chi)|^{\ell_j} \sum_{r=0}^{\ell_j} \Big( \frac{10(1+|k|)}{\ell_j}\Big)^{\ell_j-r} \frac{|k|^r}{r!}   \le \Big( \frac{12(1+|k|) |{\mathcal
P}_j(\chi)|}{\ell_j}\Big)^{\ell_j} .
\end{split}
\end{align}
  Notice that the same bound above holds for $|{\mathcal M}_j(\chi, 1-k)|$ as well. We then deduce from these estimations that when $|{\mathcal
  P}_j(\chi)| \ge \ell_j/(10(1+|k|))$, we have
\begin{align}
\label{est2}
\Big |\mathcal{N}_j(\chi,
 k)^2\mathcal{M}_j(\chi, 1-k)^{-k/(1-k)}\Big |^{2} 
& \leq \Big( \frac{12(1+|k|) |{\mathcal P}_j(\chi)|}{\ell_j}\Big)^{2(2-k/(1-k))\ell_j} \leq  |{\mathcal Q}_j(\chi, k)|^2.
\end{align}
  
  It follows \eqref{defN}, \eqref{defM}, \eqref{est1} and \eqref{est2} that we have
\begin{align*}
\sumstar_{\substack{ \chi \shortmod q }}   \Big |\mathcal{N}(\chi,
 k)^2\mathcal{M}(\chi, 1-k)^{-k/(1-k)}\Big |^{2}  \ll& \sumstar_{\substack{ \chi \shortmod q }}   \prod^R_{j=1}\Big (|{\mathcal N}_j(\chi, k)|^2 \Big( 1+ O\big(e^{-\ell_j} \big)\Big )+|{\mathcal Q}_j(\chi, k)|^2\Big ) \\
 \ll& \sumstar_{\substack{ \chi \shortmod q }}   \prod^R_{j=1}\Big (1+O( e^{-\ell_j/2})\Big ) \prod^R_{j=1}\Big (|{\mathcal N}_j(\chi, k)|^2 +|{\mathcal Q}_j(\chi, k)|^2\Big )  \\
 \ll& \sumstar_{\substack{ \chi \shortmod q }}  \prod^R_{j=1}\Big (|{\mathcal N}_j(\chi, k)|^2 +|{\mathcal Q}_j(\chi, k)|^2\Big ),
\end{align*}
  where the last estimation above follows from the observation that the sum over $e^{-\ell_j/2}$ converges. Now, a straightforward adaption of the proof of \cite[Proposition 3.5]{Gao2024-6} implies that the last expression above is $\ll \phis(q)(\log q)^{k^2}$. This completes the proof of Proposition \ref{Prop5}.

\subsection{Proof of Proposition \ref{Prop6}}

  Note that similar estimations for ${\mathcal N}_j(\chi, k)$ given in \eqref{Njboundk} and \eqref{4.2} are valid for $\mathcal{M}_j(\chi, k-1)$ as well.  We thus conclude that for each $1 \leq j \leq R$,  
\begin{align*}
\begin{split}
|{\mathcal M}_j(\chi, k-1)|^2 \leq & |{\mathcal N}_j(\chi, k-1)|^2\Big( 1+  O\Big( e^{-\ell_j} \Big )\Big )+|{\mathcal Q}_j(\chi, k)|^2 \\
\leq & \Big( 1+  O\Big( e^{-\ell_j} \Big )\Big )\Big (|{\mathcal N}_j(\chi, k-1)|^2+|{\mathcal Q}_j(\chi, k)|^2\Big ).
\end{split}
\end{align*}

  It follows that 
\begin{align}
\label{Nestmation}
\begin{split}
 \sumstar_{\substack{ \chi \shortmod q }} \big|L(\half, \chi){\mathcal M}(\chi, k-1) \big |^{2}  \leq & \sumstar_{\substack{ \chi \shortmod q }}\prod^R_{j=1}\Big (1+O( e^{-\ell_j/2})\Big )\big|L(\half, \chi)\big |^2 \prod^R_{j=1}\Big (|{\mathcal N}_j(\chi, k-1)|^2+|{\mathcal Q}_j(\chi, k)|^2\Big ) \\
 \leq & \sumstar_{\substack{ \chi \shortmod q }}\big|L(\half, \chi)\big |^2 \prod^R_{j=1}\Big (|{\mathcal N}_j(\chi, k-1)|^2+|{\mathcal Q}_j(\chi, k)|^2\Big ) \\
 \leq & \sum_{S, S^c}\sumstar_{\substack{ \chi \shortmod q }}\big|L(\half, \chi)\big |^2 \prod_{j \in S, i \in S^c}|{\mathcal N}_{j}(\chi, k-1)|^2|{\mathcal Q}_{i}(\chi, k)|^2, 
\end{split}
\end{align}
 where the sum $\sum_{S}$ is over all subsets $S$ of the set $\{ 1, \cdots, R\}$ and where we denote $S^c$ for the complement of $S$ in $\{ 1, \cdots, R\}$. 
 Now, an inspection of the proof of \cite[Proposition 3.4]{Gao2024-6} implies that for any fixed such $S$, we have
\begin{align*}
 \big|L(\half, \chi)\big |^2 \prod_{j \in S, i \in S^c}|{\mathcal N}_{j}(\chi, k-1)|^2|{\mathcal Q}_{i}(\chi, k)|^2 \ll \exp(-\sum_{i \in S^c}\frac {\ell_i}2)\phis(q)(\log q)^{ k^2 }.
\end{align*} 
 
 We deduce from this and \eqref{Nestmation} that 
\begin{align*}
\begin{split}
 \sumstar_{\substack{ \chi \shortmod q }} \big|L(\half, \chi){\mathcal M}(\chi, k-1) \big |^{2}  \ll & \prod^R_{j=1}\Big (1+O( e^{-\ell_j/2})\Big ) \phis(q)(\log q)^{ k^2 } \ll \phis(q)(\log q)^{k^2}. 
\end{split}
\end{align*}
  This establishes the desired estimation given in Proposition \ref{Prop6} and hence completes the proof. 

\vspace*{.5cm}

\noindent{\bf Acknowledgments.} The author is supported in part by NSFC grant 11871082.

\bibliography{biblio}
\bibliographystyle{amsxport}

\vspace*{.5cm}

\end{document}